\documentclass[11pt]{article}

\usepackage{enumerate}
\usepackage{amsmath}
\usepackage{amssymb,latexsym}
\usepackage{amsthm}
\usepackage{cite}
\usepackage{anysize}
\usepackage{graphicx}
\usepackage{url}
\usepackage{color}
\usepackage{multirow}

\date{}

\newcommand{\cT}{{\mathcal T}}

\newcommand{\st}{{\cal S}_t}

\newcommand{\sk}{{\cal S}_2}

\renewcommand{\proof}{\noindent{\bf Proof.}\ }
\newcommand{\Qed}{\hfill $\Box$ \medskip}

\newtheorem{theorem}{Theorem}[section]

\newtheorem{corollary}[theorem]{Corollary}
\newtheorem{definition}[theorem]{Definition}
\newtheorem{proposition}[theorem]{Proposition}

\begin{document}

\title{2-semiarcs in $\mathrm{PG}(2,q)$, $q\leq 13$}

\author{Daniele Bartoli\thanks{The research was supported by the Italian MIUR
(progetto 40\% ``Strutture Geometriche, Combinatoria e loro Applicazioni''),
and by GNSAGA.},
Giorgio Faina$^*,$ Gy\"orgy Kiss\thanks{The research was
supported by the Hungarian National
Foundation for Scientific Research, Grant No. K 81310, and by the
Slovenian-Hungarian Intergovernmental Scientific and Technological
Cooperation Project,  Grant No. T\'ET 10-1-2011-0606.}, \\
Stefano Marcugini$^*$
and Fernanda Pambianco$^*$}

\maketitle

\begin{abstract}
A $2$-semiarc is a pointset $\sk $ with the property that
the number of tangent lines to $\sk $ at each
of its points is two.
Using some theoretical results and computer aided search,
the complete classification
of $2$-semiarcs in PG$(2,q)$ is given for $q\leq 7,$
the spectrum of their sizes is determined for $q\leq 9$,
and some results about the existence are proven for $q=11$
and $q=13.$ For several sizes of $2$-semiarcs in $\mathrm{PG}(2,q)$, $q\leq 7$,
classification results have been obtained by theoretical proofs.
\end{abstract}

\section{Introduction}
\label{sec:Intro}
\indent

Ovals, $k$-arcs, and semiovals of finite projective planes are not only
interesting geometric structures, but they have important applications to
coding theory and cryptography, as well.
For details about these objects we refer the reader to
\cite{bat,HS2001, jwph1, kgy}.

Semiarcs are a  natural generalization of arcs. Let $\Pi _q$ be a projective
plane of order $q.$
A non-empty pointset ${\cal S}_t\subset \Pi _q$ is called a  {\it t-semiarc}
if for every point $P\in {\cal S}_t$
there exist exactly $t$ lines $\ell _1,\ell _2,\ldots \ell _t$ such that
${\cal S}_t \cap \ell_i = \{P\}$ for $i=1,2,\ldots ,t.$
These lines are called the \emph{tangents} to ${\cal S}_t$
at $P$. If a line $\ell $ meets $\st $ in $2,3$ or $k$ points (where $k>3$),
then $\ell $ is called a bisecant, trisecant or $k$-secant of $\st ,$
respectively.
The classical examples of semiarcs are the semiovals ($t=1$)
and the subplanes ($t=q-m,$ where $m$ is the order of the subplane).

Semiarcs are closely connected to other combinatorial structures, too.
Without the pursuit of wholeness we mention
$(r,1)$-designs and configurations.

\begin{definition}
A finite point-line
incidence structure is called \emph{linear space}
if each line contains at least two points and any
two distinct points are on exactly one line.
If there are exactly $r$ lines through each point, then the linear space
is called \emph{$(r,1)$-design}.

A \emph{$(v_r,b_k)$-configuration} is a finite point-line incidence
structure with the following properties:
\begin{itemize}
\item
There are $v$ points and $b$ lines.
\item
There are $r$ lines through each point and there are
$k$ points on each line.
\item
Two distinct lines intersect each other at most once
and two distinct points are connected by at most one line.
\end{itemize}

If $v=b$ and $r=k,$ then the configuration is called
\emph{symmetric $(v_k)$-configuration}.
\end{definition}

The following proposition gives a natural correspondence
between embeddable $(r,1)$-designs and semiarcs in finite
planes. Its proof is straightforward.

\begin{proposition}
\label{design}
If $\st $ is a $t$-semiarc in $\Pi _q,$ then the points of
$\st $ and the secants of $\st $ form a $(q+1-t,1)$-design.

If an $(r,1)$-design is embeddable to $\Pi _q,$
then its points form a $(q+1-r)$-semiarc.
\end{proposition}

$(r,1)$-designs with small $r$ were investigated by Gropp \cite{gr1, gr2}.
He constructed all $(r,1)$-designs with at most 12 points, his list
contains 974 elements, most of them are configurations. His proof is
computer assisted and he has not considered the embeddability
of these designs.

In the last years the interest and research on the fundamental problem of
determining the spectrum of the values for which there exists a given
subconfiguration
of points in $\mathrm{PG}(n,q)$ have increased considerably (see for example
\cite{ball,BDFMP2013,BDFMP2012,HS2001,KV2003,S1993,P1999,DGMP2010,DGMP2011,BDMP2011,BDFMP2013b,BFMP2013,BMP2012}).
In particular semiovals were investigated by several authors. Among others
Lisonek \cite{lis}
determined the spectrum of sizes of semiovals by exhaustive computer search
for $q\leq 9,$ $q$ odd, Bartoli \cite{bar}, Ranson and Dover
\cite{dov1, dov2},
Kiss, Marcugini and Pambianco \cite{mi1, mi2}, and Nakagawa and Suetake
\cite{nasue, sue} gave characterization theorems for semiovals in planes
of small order.

Because of the huge diversity of semiarcs, their complete classification
seems out of reach. The aim of this paper is to investigate and characterize
2-semiarcs in projective planes of order $q\leq 13.$
Throughout the paper $\Pi _q$ denotes an arbitrary projective plane of
order $q,$ while $\mathrm{PG}(2,q)$ denotes the desarguesian projective
plane over the field of $q$ elements. It is well-known,
that if $q=2,3,4,5,7$ or $8,$ then
each projective plane of order $q$ is isomorphic to $\mathrm{PG}(2,q).$

The paper is organized as follows. In Section~\ref{sec:Bounds}
we give lower and upper bounds and prove some number theoretical conditions
on the sizes of 2-semiarcs in $\Pi_q.$
Using these propositions and the results of Gropp,
in Section \ref{sec:Small}  the complete
characterization is provided for $q\leq 5.$
In Section~\ref{sec:7} we consider
the 2-semiarcs in PG$(2,7).$ A computer-free description is
given for 2-semiarcs having sizes at most 12, and a computer-assisted
proof shows that there are no 2-semiarcs in the plane with
$|\mathcal{S}_2 |\geq 13.$
Section~\ref{sec:algorithm} is devoted to the description of the
algorithm used to obtain the classification of 2-semiarcs.
Finally in Section \ref{sec:8} results about the existence of 2-semiarcs
in $\mathrm{PG}(2,q)$ for $q\in \{ 8,9,11,13\} $ are given. The computer
search is supported by the structural constraints proven in
Section \ref{sec:Bounds}.

\section{Some conditions on the sizes of 2-semiarcs}
\label{sec:Bounds}
\indent

It follows from the definition that each  $t$-semiarc in $\Pi _q$
satisfies $t\leq q+1.$ If $t$
is close to this upper bound, then we can easily classify the $t$-semiarcs.
The following proposition was proved by Csajb\'ok and Kiss \cite{csk}.

\begin{proposition}
\label{bigt}
Let $\st $ be a $t$-semiarc in $\Pi _q.$ The following properties hold:
\begin{itemize}
\item
if $t=q+1,$ then $\st $ is a single point,
\item
if $t=q,$ then $\st $ is a subset of a line, and vice versa
any subset of a
line containing at least two points is a $q$-semiarc,
\item
if $t=q-1,$ then $\st $ is a set of three non-collinear points.
\end{itemize}
\Qed
\end{proposition}

A semiarc cannot contain large collinear subsets. If
$\st $ is a $t$-semiarc in $\Pi _q,$
$\st $ is not contained in a line and it has a $k$-secant, then
$k\leq q+1-t$ obviously holds.
Semiarcs with long secants were investigated by Csajb\'ok.
He proved the following results; see \cite[Theorems 2.4 and 4.6]{cs1}.
\begin{theorem}
\label{hosszuszelo}
Let $\st $ be a $t$-semiarc in $\mathrm{PG}(2,q).$ Then the following
properties hold.
\begin{itemize}
\item
If $t<(q-1)/2,$ then $\st $ has no $(q+1-t)$-secants.
\item
If $\st $ has two $(q-t)$-secants such that the common point of these secants
is not contained in $\st $ and $\mathrm{gcd}(q,t)=\mathrm{gcd}(q-1,t-1)=1,$
then $\st $ is the union of these two $(q-t)$-secants.
\end{itemize}
\end{theorem}

Bounds on the sizes of $t$-semiarcs were also given by
Csajb\'ok and Kiss \cite{csk}. In the case
$t=2$ their result is the following.
\begin{theorem}
\label{bound}
Let $\sk$ be a $2$-semiarc in a projective plane
of order $q.$ Then
$$q\leq |\sk | \leq 1+\left\lfloor  \frac{q(1+\sqrt{8q-7})}{4}\right\rfloor .$$
\end{theorem}

\noindent
The simplest example of a $2$-semiarc of size $q$ is a $q$-arc, a set of $q$
points such that no three of them are collinear. As the following proposition
shows, there are no more examples of $2$-semiarc of size $q.$

\begin{proposition}\label{q}
Let $\sk$ be a $2$-semiarc of size $q$ in a projective plane
of order $q.$ Then $\sk$ is an arc.
\end{proposition}

{\bf Proof.}
We have to prove that no three points of $\sk $
are collinear. Suppose that the line $\ell $
is a trisecant of $\sk .$ If $P$ is a point in $\ell \cap \sk$,
 then $|\sk |=q$ implies that there are at least
$(q+1)-(q-2)=3$ tangents to $\sk $ at $P,$ contradiction.
\Qed

\begin{theorem}\label{qcon}
In $\mathrm{PG}(2,p^h)$, $p \neq 2$, there exists, up to collineations, a
unique $2$-semiarc $\sk $ of size $q=p^{h}$.
Its stabilizer group has size $hq(q-1)$.
\end{theorem}
\proof
$\sk $ is an arc of size $q=p^h.$ It is known, that
in $\mathrm{PG}(2,q)$ each $q$-arc is contained in a $(q+1)$-arc, and
if $q$ is odd, then by the Theorem of Segre,
it is contained in an irreducible conic \cite{Segre1955}.
The stabilizer of a conic is transitive on its points, hence
all the $q$-point subsets of the conic are projectively equivalent. Since
the number of conics is $q^2(q^2+q+1)(q-1)$ and each has $q+1$ subsets
of size $q,$
there are exactly $q^2(q^2+q+1)(q-1)(q+1)$ different $2$-semiarcs of size $q.$
Thus the stabilizer group has size
$\frac{|\mathrm{P\Gamma L}(3,q)|}{q^2(q^2+q+1)(q-1)(q+1)}=hq(q-1)$.
\Qed

If $\Pi _q$ contains a $2$-semiarc whose size is close to the lower bound
$q,$ then the order of the plane must satisfy some number theoretical
conditions.
\begin{proposition}
\label{q+1}
Let $\sk$ be a $2$-semiarc of size $q+1$ in a projective plane
of order $q.$ Then $q+1$ is divisible by $3.$
\end{proposition}
\proof
Let $P$ be any point of $\sk .$ The total number of lines through
$P$ is $q+1,$ and two of them are tangents to $\sk .$
The remaining $q$ points of $\sk $ are distributed among the
$q-1$ secants through $P.$ Hence there are
$q-2$ bisecants and one trisecant through $P.$
Thus each point of $\sk $ lies
on exactly one trisecant, hence $|\sk |$ is divisible by 3.
\Qed

\begin{proposition}
\label{q+2}
Let $\sk$ be a $2$-semiarc of size $q+2$ in a projective plane
of order $q.$ Then there exist integers $0\leq \alpha $ and
$0\leq \beta \neq 1$ such that $q+2=4\alpha +3\beta .$
\end{proposition}
\proof
Let $P$ be any point of $\sk .$ The total number of lines through
$P$ is $q+1,$ two of them are tangents to $\sk .$
The remaining $q+1$ points of $\sk $ are distributed among the
$q-1$ secants through $P.$ Hence there are
either two trisecants and $q-3$ bisecants,
or one $4$-secant and $q-2$ bisecants through $P.$
Thus each point lies on either two trisecants or one $4$-secant.
Let $\mathcal{T}_{3}$ be the set of points lying on two trisecants. Then
it is a configuration $(v_2,k_3)$, where $v=|\mathcal{T}_{3}|$ and,
by  \cite[Theorem 3.1]{gr}, $|\mathcal{T}_{3}|=3\beta$, with $\beta\neq 1$.
Let $\mathcal{T}_{4}$
be the set of points lying on one $4$-secant, then $|\mathcal{T}_{4}|=4\alpha$.
Then $q+2=4\alpha +3\beta$, with $\alpha, \beta \geq 0$ and $\beta \neq 1$.
\Qed

\section{The small planes}
\label{sec:Small}
\indent

The classification of
$2$-semiarcs in the cases
$q=2$ and $q=3$ follows from Proposition \ref{bigt}.

\begin{theorem} $ $
\begin{itemize}
\item
In $\mathrm{PG}(2,2)$ each $2$-semiarc $\sk$ consists of two or three
collinear points.
\item
In $\mathrm{PG}(2,3)$ each $2$-semiarc $\sk$ is a set of three non-collinear
points.
\end{itemize}
\end{theorem}

If $q=4,$ then 2-semiarcs correspond to $(3,1)$-designs by Proposition \ref{design}.
Gropp \cite[Table 1]{gr1} proved that there are three such designs, they consist of
4, 6 and 7 points, respectively. He also gave a detailed combinatorial description
of these objets. We show that each of these designs is embeddable into
$\mathrm{PG}(2,4).$

\begin{theorem}
In $\mathrm{PG}(2,4)$ there are three projectively non-equivalent
$2$-semiarcs.
\begin{itemize}
\item
$|\sk |=4,$ four points in general position.
\item
$|\sk |=6,$ the vertices of a complete quadrilateral.
\item
$|\sk |=7,$ the points of a subplane of order 2.
\end{itemize}
\end{theorem}

{\bf Proof.}
It is easy to verify (without applying Gropp's results), that there are only
three possible sizes of a 2-semiarc.
Theorem \ref{bound} gives $4\leq |\sk |\leq 7.$ From Proposition \ref{q+1}
we get $|\sk |\neq 5,$ because $q+1=5$ is not divisible by 3.
Hence $|\sk |\in \{ 4,6,7\} .$

The case $|\sk |=4$ follows from Proposition \ref{q}.
The combinatorial description of Gropp gives that if $|\sk |=6,$ then
there are two trisecants and one bisecant through each point, hence the
design corresponds to the six vertices of a complete quadrilateral and
it is obviously embeddable into $\mathrm{PG}(2,4).$
If $|\sk |=7,$ then according to Gropp, the design is a $(7_3)$-configuration.
In other words this is the Fano plane $\mathrm{PG}(2,2),$ which is embeddable to
$\mathrm{PG}(2,4).$
\Qed

\begin{table}[h]
\begin{center}
\begin{tabular}{|c|c|c|c|c|c|}
\hline
$|\mathcal{S}_2|$&$x_0$&$x_1$&$x_2$&$x_3$&$G$\\
\hline
\hline
$4$ &
7 &
8 &
6 &
0 &
$\mathbb{Z}_2 \times \mathrm{S}_4$\\
\hline
$6$ &
2 &
12 &
3 &
4 &
$\mathbb{Z}_2 \times \mathrm{S}_4$\\
\hline
$7$ &
0 &
14 &
0 &
7 &
$\mathrm{PSL}(3,2)\times \mathbb{Z}_2$ \\
\hline
\end{tabular}
\end{center}
\begin{center}
\caption{$2$-semiarcs in $\mathrm{PG}(2,4)$}\label{table:PG2_4_2}
\end{center}
\end{table}

Table \ref{table:PG2_4_2} contains the non-equivalent $2$-semiarcs in
$\mathrm{PG}(2,4)$, the number of their $i$-secants, $x_i$,
and the description
of the stabilizer groups in $\mathrm{P\Gamma L} (3,4)$.


If $q=5,$ then 2-semiarcs correspond to $(4,1)$-designs by Proposition \ref{design}.
Gropp \cite[Table 1]{gr1} proved that there are eight such designs with at most 12
points. We show that only three of them are embeddable into
$\mathrm{PG}(2,5).$

\begin{theorem}
In $\mathrm{PG}(2,5)$ there are three projectively non-equivalent
$2$-semiarcs.
\begin{itemize}
\item
$|\sk |=5,$ five points of a conic.
\item
$|\sk |=6,$ the union of two trisecants.
\item
$|\sk |=9,$ the projective triangle.
\end{itemize}
\end{theorem}

{\bf Proof.}
It is easy to see that there are only
four possible sizes of a 2-semiarc.
Theorem \ref{bound} gives $5\leq |\sk |\leq 9.$ From Proposition \ref{q+2}
we get $|\sk |\neq 7,$ because $q+2=7$ cannot be written as $4\alpha +3\beta $
with $\beta \neq 1.$

First we prove that $|\sk |\neq 8.$ Suppose to the contrary
that $\sk $ is a $2$-semiarc with 8 points. Gropp proved that
there is only one $(4,1)$-design with eight points, the
symmetric $(8_3)$-configuration (also called M\"obius-Kantor
configuration). But it was proven by
Abdul-Elah, Al-Dhahir and Jungnickel \cite{aaj} that this configuration
cannot be embedded into $\mathrm{PG}(2,5).$
Hence $|\sk |\in \{ 5,6,9\} .$

The case $|\sk |=5$ follows from Proposition \ref{q}.

In the case $|\sk |=6$,
if $P\in \sk $ is a point, then there
are $6-2=4$ non-tangents through $P,$ hence
$\sk $ has no 4-secants. Let
$a$ be the number of trisecants, and $b$ be the number
of bisecants  through $P.$ Then we get $a+b=4$ and $2a+b=5,$ hence $a=1$
and $b=3.$ So $\sk $ is the union of two trisecants, $\ell _1$ and $\ell _2.$
This is the second case of Theorem \ref{hosszuszelo}.

Finally consider the case $|\sk |=9.$ Gropp proved that
there are two $(4,1)$-designs with nine points. One of them
is the affine plane of order 3. But $\mathrm{AG}(2,3)$ cannot be embedded
into $\mathrm{PG}(2,q)$ if $q\equiv 2$ (mod 3) (see e.g. \cite{bk}).

The points of the other $(4,1)$-design are of two types: (i) the vertices of
a triangle $\cT ,$ (ii) the points on exactly one side of $\cT ,$ two points on
each side.
If a point is of type (i), then it is on two 4-secants and
on two bisecants; if a point is of type (ii), then it is on one 4-secant and hence
on two trisecants and on one bisecant.
Hence $\sk $ has three 4-secants, $6\cdot 2/3=4$ trisecants and $(3\cdot 2+6\cdot 1)/2=6$
bisecants. $\sk $ also has $9\cdot 2=18$ tangents, so $\sk $ is a blocking set
because $3+4+6+18=31$ equals to the total number of lines in $\mathrm{PG}(2,5).$
This blocking set has cardinality $3(q+1)/2,$ hence by a theorem of Lov\'asz and Schrijver
\cite{ls} it is a projective triangle.

A possible embedding into
$\mathrm{PG}(2,q)$ is the following. The vertices of $\cT $:
$\{(1:0:0),(0:1:0),(0:0:1)\},$ the points on the sides of $\cT $:
$\{ (1:1:0),(4:1:0),
(1:0:1),(4:0:1),(0:1:1),(0:4:1) \} .$
\Qed

\begin{table}[h]
\begin{center}
\begin{tabular}{|c|c|c|c|c|c|c|}
\hline
$|\mathcal{S}_2|$&$x_0$&$x_1$&$x_2$&$x_3$&$x_{4}$&$G$\\
\hline
\hline
$5$ &
11 &
10 &
10 &
0&
0&
$\mathbb{Z}_{5}\rtimes \mathbb{Z}_{4}$\\
\hline
$6$ &
8 &
12 &
9 &
2 &
0&
$\mathrm{D}_4$ \\
\hline
$9$ &
0 &
18 &
6 &
4 &
3 &
$\mathrm{S}_{4}$\\
\hline
\end{tabular}
\end{center}
\begin{center}
\caption{$2$-semiarcs in $\mathrm{PG}(2,5)$}\label{table:PG2_5_2}
\end{center}
\end{table}

Table \ref{table:PG2_5_2} contains the non-equivalent $2$-semiarcs
in $\mathrm{PG}(2,5)$,
the number of their $i$-secants, $x_i$,
and the description of the stabilizer groups in $\mathrm{PGL} (3,5)$.

\section{$2$-semiarcs in $\mathrm{PG}(2,7)$}
\label{sec:7}
\indent

The number of $(6,1)$-designs with at most
12 points is 47. Instead of considering the list
of Gropp \cite{gr1}, we give a geometric characterization
of the embeddable designs and we prove that there are
25 non-equivalent 2-semiarcs in $\mathrm{PG}(2,7).$
First consider the long secants of the semiarcs.
If $q=7$ and $t=2,$ then Theorem \ref{hosszuszelo} gives
the following corollary.

\begin{corollary}
\label{56}
Let $\sk $ be a $2$-semiarc in $\mathrm{PG}(2,7).$ Then
$\sk $ has no $6$-secants.
If $\sk $ has two $5$-secants such that the common point of these secants
is not contained in $\sk ,$
then $\sk $ is the union of these two $5$-secants.
\end{corollary}

\noindent
If the common point of the long secants belongs to $\sk ,$ then the size of
the semiarc cannot be small.

\begin{proposition}
\label{biggerthan12}
Let $\sk $ be a $2$-semiarc in $\mathrm{PG}(2,7).$
If $\sk $ has two $5$-secants such that the common point of these secants
is contained in $\sk ,$ then $|\sk |>12.$
\end{proposition}

{\bf Proof.}
Let $\ell _1$ and $\ell _2$ be the $5$-secants and let
$P\in \ell _1\cap \ell _2.$
Then $P\in \sk $ implies that there are six secants of $\sk $
through $P.$ Hence $\sk \setminus (\ell _1\cup \ell _2)$ must
contain at least four points. So $|\st |\geq 9+4=13$ holds.
\Qed

\begin{center}
\begin{table}[t]
\caption{$2$-semiarcs in $\mathrm{PG}(2,7)$}\label{table:PG2_7_2}
\begin{center}
\begin{tabular}{|c|c|c|c|c|c|c|c|c|c|}
\hline
$|\mathcal{S}_2|$&$\mathcal{S}_2$&$x_0$&$x_1$&$x_2$&$x_3$&$x_4$&$x_5$&$G$\\
\hline
\hline
$7$&
\tabcolsep=0.85 mm
{\tiny
\begin{tabular}{ccccccc}
 1& 1& 1& 0& 0& 1& 1\\
 1& 2& 0& 1& 0& 3& 6\\
 1& 3& 0& 0& 1& 2& 4\\
\end{tabular}}&
22 &
14 &
21 &
0&
0&
0&
$G_{42}$\\
\hline
\hline
$9$&{\tiny
\tabcolsep=0.85 mm
\begin{tabular}{ccccccccc}
 0& 1& 0& 1& 0& 1& 1& 1& 1\\
 1& 0& 1& 3& 0& 6& 1& 5& 1\\
 5& 0& 0& 2& 1& 3& 1& 2& 5\\
\end{tabular}}&
15 &
18 &
18 &
6 &
0&
0&
$\mathbb{Z}_{2}$\\
\hline
$9$&{\tiny
\tabcolsep=0.85 mm
\begin{tabular}{ccccccccc}
 0& 1& 0& 1& 0& 1& 1& 1& 1\\
 1& 0& 1& 3& 0& 4& 0& 6& 1\\
 5& 0& 0& 2& 1& 1& 4& 5& 1\\
\end{tabular}}&
15 &
18 &
18 &
6 &
0&
0&
$\mathbb{Z}_{2}$\\
\hline
$9$&{\tiny
\tabcolsep=0.85 mm
\begin{tabular}{ccccccccc}
 0& 1& 0& 1& 0& 1& 1& 1& 1\\
 1& 0& 1& 3& 0& 6& 4& 6& 1\\
 5& 0& 0& 2& 1& 3& 0& 5& 1\\
\end{tabular}}&
15 &
18 &
18 &
6 &
0&
0&
$\mathbb{Z}_{3}$\\
\hline
$9$&{\tiny
\tabcolsep=0.85 mm
\begin{tabular}{ccccccccc}
 0& 1& 0& 1& 0& 1& 1& 1& 1\\
 1& 0& 1& 3& 0& 4& 3& 6& 1\\
 5& 0& 0& 2& 1& 1& 3& 5& 1\\
\end{tabular}}&
15 &
18 &
18 &
6 &
0&
0&
$\mathbb{Z}_{3}$\\
\hline
$9$&{\tiny
\tabcolsep=0.85 mm
\begin{tabular}{ccccccccc}
 1& 1& 1& 0& 0& 1& 0& 1& 1\\
 1& 3& 0& 1& 0& 3& 1& 5& 2\\
 1& 5& 0& 0& 1& 2& 5& 1& 6\\
\end{tabular}}&
15 &
18 &
18 &
6 &
0&
0&
$\mathbb{Z}_{6}$\\
\hline
$9$&
{\tiny
\tabcolsep=0.85 mm
\begin{tabular}{ccccccccc}
 1& 0& 1& 0& 1& 0& 1& 1& 1\\
 4& 1& 0& 1& 3& 0& 6& 1& 1\\
 3& 5& 0& 0& 2& 1& 3& 1& 5\\
\end{tabular}}&
15 &
18 &
18 &
6 &
0&
0&
$\mathrm{S}_{3}$\\
\hline
\hline
$10$&{\tiny
\tabcolsep=0.85 mm
\begin{tabular}{cccccccccc}
 1& 1& 0& 0& 1& 1& 0& 1& 1& 1\\
 1& 0& 1& 0& 3& 5& 1& 6& 1& 6\\
 1& 0& 0& 1& 2& 0& 5& 6& 2& 5\\
\end{tabular}}&
12 &
20 &
15 &
10 &
0 &
0 &
$\mathbb{Z}_{1}$\\
\hline
$10$&{\tiny
\tabcolsep=0.85 mm
\begin{tabular}{cccccccccc}
 1& 1& 1& 0& 0& 1& 1& 0& 1& 1\\
 1& 3& 0& 1& 0& 3& 5& 1& 2& 1\\
 1& 5& 0& 0& 1& 2& 0& 5& 1& 3\\
\end{tabular}}&
12 &
20 &
15 &
10 &
0 &
0 &
$\mathbb{Z}_{1}$\\
\hline
$10$&{\tiny
\tabcolsep=0.85 mm
\begin{tabular}{cccccccccc}
 1& 1& 0& 0& 1& 1& 0& 1& 1& 1\\
 1& 0& 1& 0& 3& 3& 1& 0& 4& 6\\
 1& 0& 0& 1& 2& 6& 5& 1& 2& 5\\
\end{tabular}}&
12 &
20 &
15 &
10 &
0 &
0 &
$\mathbb{Z}_{1}$\\
\hline
$10$&{\tiny
\tabcolsep=0.85 mm
\begin{tabular}{cccccccccc}
 1& 1& 1& 0& 0& 1& 1& 1& 0& 1\\
 1& 3& 0& 1& 0& 3& 0& 5& 1& 2\\
 1& 5& 0& 0& 1& 2& 6& 0& 5& 1\\
\end{tabular}}&
12 &
20 &
15 &
10 &
0 &
0 &
$\mathbb{Z}_{2}$\\
\hline
$10$&{\tiny
\tabcolsep=0.85 mm
\begin{tabular}{cccccccccc}
 1& 1& 1& 0& 0& 1& 1& 1& 0& 1\\
 1& 4& 0& 1& 0& 3& 0& 1& 1& 6\\
 1& 1& 0& 0& 1& 2& 3& 0& 5& 5\\
\end{tabular}}&
12 &
20 &
15 &
10 &
0 &
0 &
$\mathbb{Z}_{2}$\\
\hline
$10$&{\tiny
\tabcolsep=0.85 mm
\begin{tabular}{cccccccccc}
 0& 1& 0& 1& 0& 1& 1& 1& 1& 1\\
 1& 0& 1& 3& 0& 0& 2& 3& 1& 5\\
 5& 0& 0& 2& 1& 6& 6& 5& 1& 5\\
\end{tabular}}&
12 &
20 &
15 &
10 &
0 &
0 &
$\mathbb{Z}_{2}$\\
\hline
$10$&{\tiny
\tabcolsep=0.85 mm
\begin{tabular}{cccccccccc}
 1& 1& 1& 0& 0& 1& 1& 0& 1& 1\\
 1& 6& 0& 1& 0& 3& 1& 1& 4& 6\\
 1& 3& 0& 0& 1& 2& 0& 5& 2& 5\\
\end{tabular}}&
12 &
20 &
15 &
10 &
0 &
0 &
$\mathbb{Z}_{2}$\\
\hline
$10$&{\tiny
\tabcolsep=0.85 mm
\begin{tabular}{cccccccccc}
 1& 1& 0& 0& 1& 1& 1& 0& 1& 1\\
 1& 0& 1& 0& 2& 3& 3& 1& 0& 6\\
 1& 0& 0& 1& 2& 2& 6& 5& 1& 5\\
\end{tabular}}&
12 &
20 &
15 &
10 &
0 &
0 &
$\mathbb{Z}_{3}$\\
\hline
$10$&{\tiny
\tabcolsep=0.85 mm
\begin{tabular}{cccccccccc}
 0& 1& 1& 0& 1& 0& 1& 1& 1& 1\\
 1& 0& 0& 1& 3& 0& 4& 6& 1& 1\\
 5& 0& 1& 0& 2& 1& 2& 5& 1& 3\\
\end{tabular}}&
12 &
20 &
15 &
10 &
0 &
0 &
$\mathbb{Z}_{4}$\\
\hline
$10$&{\tiny
\tabcolsep=0.85 mm
\begin{tabular}{cccccccccc}
 0& 1& 0& 0& 1& 0& 1& 1& 0& 1\\
 1& 0& 1& 0& 3& 1& 1& 5& 1& 2\\
 5& 0& 0& 1& 3& 3& 1& 5& 2& 2\\
\end{tabular}}&
10 &
20 &
25 &
0 &
0 &
2 &
$\mathrm{D}_{4}$\\
\hline
$10$&{\tiny
\tabcolsep=0.85 mm
\begin{tabular}{cccccccccc}
 1& 1& 0& 0& 1& 1& 1& 0& 1& 1\\
 1& 0& 1& 0& 3& 0& 5& 1& 5& 2\\
 1& 0& 0& 1& 2& 6& 4& 5& 1& 6\\
\end{tabular}}&
12 &
20 &
15 &
10 &
0 &
0 &
$\mathrm{D}_{6}$\\
\hline
$10$&{\tiny
\tabcolsep=0.85 mm
\begin{tabular}{cccccccccc}
 0& 1& 1& 0& 1& 0& 1& 1& 1& 1\\
 1& 0& 6& 1& 3& 0& 3& 0& 1& 5\\
 5& 0& 4& 0& 2& 1& 3& 4& 1& 2\\
\end{tabular}}&
12 &
20 &
15 &
10 &
0 &
0 &
$\mathcal{S}_{4}$\\
\hline
\hline
$11$&{\tiny
\tabcolsep=0.85 mm
\begin{tabular}{ccccccccccc}
 0& 1& 0& 1& 0& 1& 0& 1& 1& 1& 1\\
 1& 0& 1& 3& 0& 1& 1& 6& 1& 5& 1\\
 5& 0& 0& 2& 1& 6& 3& 5& 1& 0& 5\\
\end{tabular}}&
8 &
22 &
19 &
4 &
4 &
0&
$\mathbb{Z}_{1}$\\
\hline
$11$&{\tiny
\tabcolsep=0.85 mm
\begin{tabular}{ccccccccccc}
 0& 1& 0& 1& 1& 0& 1& 1& 1& 1& 1\\
 1& 0& 1& 6& 3& 0& 1& 4& 6& 0& 1\\
 5& 0& 0& 6& 2& 1& 2& 2& 5& 5& 1\\
\end{tabular}}&
9 &
22 &
13 &
12 &
1 &
0&
$\mathbb{Z}_{1}$\\
\hline
$11$&{\tiny
\tabcolsep=0.85 mm
\begin{tabular}{ccccccccccc}
 0& 1& 1& 0& 1& 0& 1& 1& 1& 1& 1\\
 1& 0& 0& 1& 3& 0& 1& 6& 1& 6& 5\\
 5& 0& 1& 0& 2& 1& 6& 5& 1& 2& 5\\
\end{tabular}}&
9 &
22 &
13 &
12 &
1 &
0&
$\mathbb{Z}_{1}$\\
\hline
\hline
$12$&{\tiny
\tabcolsep=0.85 mm
\begin{tabular}{cccccccccccc}
 1& 0& 1& 0& 1& 0& 1& 1& 1& 0& 1& 1\\
 0& 1& 1& 0& 3& 1& 1& 6& 2& 1& 0& 2\\
 0& 0& 3& 1& 2& 3& 1& 5& 2& 5& 3& 1\\
\end{tabular}}&
6 &
24 &
12 &
12 &
3 &
0&
$\mathbb{Z}_{1}$\\
\hline
$12$&{\tiny
\tabcolsep=0.85 mm
\begin{tabular}{cccccccccccc}
 0& 1& 0& 1& 0& 1& 1& 1& 1& 1& 1& 1\\
 1& 0& 1& 3& 0& 4& 3& 6& 1& 6& 1& 0\\
 5& 0& 0& 2& 1& 0& 0& 5& 1& 2& 3& 3\\
\end{tabular}}&
6 &
24 &
12 &
12 &
3 &
0&
$\mathbb{Z}_{3}$\\
\hline
$12$&{\tiny
\tabcolsep=0.85 mm
\begin{tabular}{cccccccccccc}
 0& 1& 0& 1& 0& 1& 1& 1& 0& 1& 1& 1\\
 1& 0& 1& 3& 0& 3& 6& 1& 1& 6& 1& 5\\
 5& 0& 0& 2& 1& 0& 5& 1& 1& 2& 3& 0\\
\end{tabular}}&
6 &
24 &
12 &
12 &
3 &
0&
$\mathbb{Z}_{3}$\\
\hline
\end{tabular}
\end{center}
\end{table}
\end{center}

\begin{theorem}
In $\mathrm{PG}(2,7)$ there are nine combinatorially
non-equivalent $2$-semiarcs (there are projectively non-equivalent
subclasses in some combinatorial classes).
\begin{itemize}
\item
$|\sk |=7,$ seven points of a conic.
\item
$|\sk |=9,$ there are two types,
\begin{enumerate}
\item
nine vertices of a $3\times 3$ grid,
\item
the six vertices of two triangles $\cT _1$ and $\cT _2,$
and the three points of
intersections of the corresponding sides of $\cT _1$ and $\cT _2.$
\end{enumerate}
\item
$|\sk |=10,$ there are two types,
\begin{enumerate}
\item
the union of two 5-secants,
\item
the points of a $10_3$ configuration.
\end{enumerate}
\item
$|\sk |=11,$ then the semiarc has no 5-secant. There are two types,
\begin{enumerate}
\item
four 4-secants and four trisecants,
\item
one 4-secant and ten trisecants.
\end{enumerate}
\item
$|\sk |=12,$ then it has three 4-secants and these lines form a
triangle $\mathcal T$. There are two types,
\begin{enumerate}
\item
two vertices of $\mathcal T$ belong to $\sk ,$
\item
three vertices of $\mathcal T$ belong to $\sk .$
\end{enumerate}
\end{itemize}
\end{theorem}

{\bf Proof.}
Theorem \ref{bound} gives $7\leq |\sk |\leq 15.$
Let $s$ be the number of points of ${\cal S}_2,$ let
${\cal L}=\{ \ell_1,\ell_2,\ldots ,\ell_{57}\} $
be the set of lines of $\mathrm{PG}(2,7)$ and
let $c_i=|{\cal S}_2\cap \ell _i|$ for $i=1,2,\ldots ,57.$
If we count in two
different ways the number of
incident point-line pairs $(P,\ell _j)$ where $\ell _j \in {\cal L}$ and $P\in {\cal S}_2,$
and the ordered triples $(P_1,P_2, \ell _j)$ where $\ell _j\in {\cal L}$ and the
distinct points $P_1$
and $P_2$ are in ${\cal S}_2\cap \ell _j,$
then we get
$$\sum _{i=1}^{57}c_i = 8s \quad {\rm and} \quad
\sum _{i=1}^{57}c_i(c_i-1)=s(s-1).$$
Hence
$$\sum _{i=1}^{57}c_i^2 = s^2+7s.$$
We may assume without loss of generality that the lines
$\ell _{58-2s}, \ell _{59-2s},\ldots ,\ell _{57}$
are the tangents to ${\cal S}_2,$ for these lines $c_i=1.$ If we subtract
these values, then we get
\begin{equation}
\label{becsles}
\sum _{i=1}^{57-2s}c_i = 6s \quad {\rm and} \quad
\sum _{i=1}^{57-2s}c_i^2 = s^2+5s.
\end{equation}

It follows from Corollary \ref{56} that if $k\geq 6,$ then $\sk $ has no $k$-secant.
Let $x_i$ be the number of $i$-secants of ${\cal S}_2$ for $i=0,1,\ldots ,5.$
Then
$$\sum _{i=1}^{57-2s}(c_i-2)(c_i-3)=6x_0+2x_4+6x_5 \quad {\rm and} \quad
\sum _{i=1}^{57-2s}(c_i-3)(c_i-4)=12x_0+2x_2+2x_5.$$
On the other hand, Equations (\ref{becsles}) give
$$\sum _{i=1}^{57-2s}(c_i-2)(c_i-3)=\sum _{i=1}^{57-2s}(c_i^2-5c_i+6)=s^2+5s-5\cdot 6s+6(57-2s)=s^2-37s+342$$
and
$$\sum _{i=1}^{57-2s}(c_i-3)(c_i-4)=\sum _{i=1}^{57-2s}(c_i^2-7c_i+12)=s^2+5s-7\cdot 6s+12(57-2s)=s^2-61s+684.$$
Hence
\begin{equation}
\label{korlat}
6x_0+2x_4+6x_5=s^2-37s+342 \quad {\rm and} \quad
12x_0+2x_2+2x_5=s^2-61s+684.
\end{equation}

First we prove the non-existence parts of the theorem.
From Proposition \ref{q+1}
we get $|\sk |\neq 8,$ because $q+1=8$ is not divisible by 3.

Suppose, that $s=15.$ Then Equations (\ref{becsles}) give
$\sum _{i=1}^{27}c_i = 90$ and
$\sum _{i=1}^{27}c_i^2 =300.$
Applying the inequality between the arithmetic and quadratic means we get
$$ \frac{90}{27}=\frac{\sum _{i=1}^{27}c_i}{27}
\le \sqrt{\frac{\sum _{i=1}^{27}c_i^2}{27}}
= \sqrt{\frac{300}{27}}=\frac{10}{3}.$$
Thus equality holds, hence $c_1=c_2=\ldots =c_{27}.$ But
$90/27$ is not an integer, contradiction.

Now suppose, that $s=14.$ Then Equations (\ref{korlat}) give
$$3x_0+x_4+3x_5=10 \quad {\rm and} \quad 6x_0+x_2+x_5=13.$$
Elementary counting shows that there are only nine possibilities
for the numbers $x_0,x_1,\ldots ,x_5.$ These are the following.

\smallskip
\begin{center}
\begin{tabular}{|c|c|c|c|c|c|}
\hline
$x_0$&$x_1$&$x_2$&$x_3$&$x_4$&$x_5$ \\
\hline
\hline
2&28&1&22&4&0 \\
\hline
1&28&7&14&7&0 \\
\hline
0&28&13&6&10&0 \\
\hline
2&28&0&25&1&1 \\
\hline
1&28&6&17&4&1 \\
\hline
0&28&12&9&7& 1 \\
\hline
1&28&5&20&1&2 \\
\hline
0&28&11&12&4&2 \\
\hline
0&28&10&15&1& 3 \\
\hline
\end{tabular}
\end{center}

Now suppose, that $s=13.$ Then Equations (\ref{korlat}) give
$$3x_0+x_4+3x_5=15 \quad {\rm and} \quad 6x_0+x_2+x_5=30.$$
Elementary counting shows that there are only twelve possibilities
for the numbers $x_0,x_1,\ldots ,x_5.$ These are as follows.

\smallskip
\begin{center}
\begin{tabular}{|c|c|c|c|c|c|}
\hline
$x_0$&$x_1$&$x_2$&$x_3$&$x_4$&$x_5$ \\
\hline
\hline
5&26&0&26&0&0 \\
\hline
4&26&6&18&3&0 \\
\hline
3&26&12&10&6&0 \\
\hline
2&26&18&2&9&0 \\
\hline
4&26&5&21&0&1 \\
\hline
3&26&11&13&3& 1 \\
\hline
2&26&17&5&6& 1 \\
\hline
3&26&10&16&0& 2 \\
\hline
2&26&16&8&3& 2 \\
\hline
1&26&22&0&6& 2 \\
\hline
2&26&15&11&0& 3 \\
\hline
1&26&21&3&3& 3 \\
\hline
\end{tabular}
\end{center}

\noindent
In these cases an exhaustive computer search shows that there are no
$2$-semiarcs of sizes 14 and 13 in $\mathrm{PG}(2,7).$

Now consider the existence parts.
The case $|\sk |=7$ follows from Proposition \ref{q}.

\smallskip
If $|\sk |=9$ then we can apply Proposition \ref{q+2}.
As $9=4\alpha +3\beta $ implies $\alpha =0$ and $\beta =3,$
we get that there is no 4-secant of $\sk $ and there are
two trisecants through each point of $\sk .$ Hence the total number of trisecants is $9\times 2/3=6.$
There are two possibilities.

(i) There do not exist three trisecants such that they form a triangle
whose three vertices are in $\sk .$
Then the points of $\sk $ are the nine vertices of a $3\times 3$ grid,
whose six lines are the trisecants of $\sk .$ An example for this case is the following.
The points of $\sk $ are the points of intersections of three horizontal
and three vertical lines. Their cartesian coordinates are the following:
$(0,0),\, (1,0),\, (3,0),$
$(0,1),\, (1,1),\, (3,1),$ $(0,4),\, (1,4)$ and $(3,4).$

The grid has two triples of lines. There are two possibilities in each triples:
the lines either form a triangle or they belong to a pencil.
Hence there are projectively non-isomorphic examples of this combinatorial type (see Table \ref{table:PG2_7_2}).

(ii) There exist three trisecants such that they form a triangle
$\cT _1$ whose three vertices, say $P_1,P_2$ and $P_3$ are in $\sk .$
In this case $\sk $ contains three points, say $Q_1,Q_2$ and $Q_3$
from the sides of $\cT _1,$ and three more points, say $R_1,R_2$ and $R_3.$
Consider the three other trisecants of $\sk .$ If $Q_iQ_j$ were a trisecant,
then it ought to contain exactly one point from the set $\{ R_1,R_2,R_3\} ,$
hence both of the remaining two trisecants would pass on the other two $R_i,$
contradiction. So each of the remaining three trisecants contains one point
of the set $\{ Q_1,Q_2,Q_3\} ,$ hence two points from the
set $\{ R_1,R_2,R_3\} .$
So the points $R_1,R_2$ and $R_3$ form a triangle $\cT _2.$
An example for this case is the following.
The homogeneous coordinates of the vertices of $\cT _1$ are $(0:0:1),\,(0:1:0)$ and $(1:0:0),$ the coordinates of
the vertices of $\cT _2$ are $(2:3:1),\, (3:4:1)$ and $(5:5:1).$ The points of
intersections of the corresponding sides are $(1:4:0),\, (0:1:1)$ and $(1:0:1).$

There are projectively non-isomorphic examples of this combinatorial
type, too (see Table \ref{table:PG2_7_2}).

\smallskip
If $|\sk |=10,$ then first we consider the largest collinear subset
of $\sk .$ Because of Theorem \ref{hosszuszelo}
its cardinality is at most $q-2=5.$ If $\sk $ has a $5$-secant,
then Csajb\'ok, H\'eger and Kiss \cite[Proposition 2.3]{cshk}
proved that $\sk $ is the union of two 5-secants.

If $\sk $ has no 5-secants, then
the points of $\sk $ can be partitioned into
two subsets. Let ${\mathcal A}\subset \sk $ be the set of points
belonging to three trisecants of $\sk $ and let ${\mathcal B}\subset \sk $
be the set of points belonging to one trisecant  and one 4-secant of $\sk .$
If $|{\mathcal A}|=a$ and $|{\mathcal B}|=b,$ then the total number of
trisecants of $\sk $ is $(3a+b)/3,$ hence $3|b.$ Thus if $b>0,$ then
$b\geq 3,$ and no point of $\sk $ lies on more than one 4-secant.
Hence $b>0$ implies $|\sk |\geq 3\times 4=12,$ contradiction.
So $\sk $ has no 4-secant, hence it is a $(10_3)$-configuration.
An example for this case is the Desargues configuration.

It is known that there are ten projectively non-isomorphic $(10_3)$-configurations
\cite {gr}. The embeddability of these configurations were
investigated by Glynn \cite{glynn}, who proved that one of them is
not embeddable into any pappian plane. It is also known, that the other
nine can be embedded into the classical euclidean plane \cite{bokstu}.
Our exhaustive computer search shows that these nine can also be
embedded into $\mathrm{PG}(2,7)$.


\smallskip
If $|\sk |=11,$ then it is a $2$-semiarc with $q+4$
points. For each point $P\in \sk $ there are $q-1$
secants through $P,$ thus
$q+3$ points of $\sk $ are distributed
among the secants through $P.$
It follows from Corollary \ref{56} that
$\sk $ has no 6-secant.
Thus the points of $\sk $ can be partitioned into
four subsets. Let ${\mathcal A}\subset \sk $ be the set of points
belonging to four trisecants of $\sk ,$ let ${\mathcal B}\subset \sk $
be the set of points belonging to two trisecants and
one 4-secant of $\sk ,$
let ${\mathcal C}\subset \sk $
be the set of points belonging to two 4-secants of $\sk $ and finally
let ${\mathcal D}\subset \sk $
be the set of points belonging to one trisecant and one $5$-secant of $\sk .$

First we prove that ${\mathcal D}=\emptyset .$
Let $|{\mathcal A}|=a,$ $|{\mathcal B}|=b,$
$|{\mathcal C}|=c$ and $|{\mathcal D}|=d.$
Let $s$ be the number of $5$-secants.
Then Corollary \ref{56} and Proposition \ref{biggerthan12} imply that
$s\leq 1.$ Suppose that $s=1.$ Then we show that $c=0$
also holds. The 4-secants cannot meet the
5-secant in a point of $\sk$ and the union of two intersecting 4-secants
contains 7 points, so if
$c\neq 0,$ then $\sk $ contains at least $5+7>11$ points, contradiction.
So $s=1$ implies $a+b=6.$
The number of
the 4-secants of $\sk $ is $b/4,$ hence $4|b.$ There are only
two possibilities, either $b=0$ or $b=4.$
In the first case $a=6,$ in the second $a=2.$  The number of
the trisecants of $\sk $ is $(5+4a+2b)/3$.
If $b=0$, then this number is $5+4\cdot 6=29$ and it is not divisible by 3, contradiction.
If $b=4$ then  $a=2,$ and $\sk $ has one 5-secant, $\ell _5 ,$ one 4-secant, $\ell _4$ and seven trisecants.
Let $\sk \setminus (\ell _5 \cup \ell _4 )=\{ P,R\}.$
Then there are four trisecants through both $P$ and $R,$ hence the
line $PR$ is a trisecant. Each of the other $2\times 3=6$ trisecants through
$P$ or $R$ must contain one point of $\ell _5,$ but there exists a unique
trisecant  at each point of $\ell _5.$ This contradiction proves $d=0.$

If $d=0,$ then $a+b+c=11.$
The number of
the trisecants of $\sk $ is $(4a+2b)/3,$ hence
$$b\equiv a\quad ({\rm mod} \, 3).$$
The number of
the 4-secants of $\sk $ is $(b+2c)/4=(22-2a-b)/4,$ hence
$$b\equiv 2a+2 \quad ({\rm mod} \, 4).$$
Thus the Chinese Remainder Theorem gives
$$b\equiv 10a+6 \quad ({\rm mod} \, 12).$$
We know that $0\leq a,b\leq 11,$ hence if $a$ is
given, then this congruence uniquely determines
$b,$ and also $0\leq c=11-a-b.$ We have the following
possibilities.

$$
\begin{array}{|r||c|c|c|c|c|c|c|c|c|c|c|c|}
\hline
a & 0 & 1 & 2 & 3 & 4 & 5 & 6 & 7 & 8 & 9 & 10 & 11 \\
\hline
b & 6 & 4 & 2 & 0 & 10 & 8 & 6 & 4 & 4 & 0 & 10 & 8 \\
\hline
c & 5 & 6 & 7 & 8 & - & - & - & 0 & - & 2 & - & -  \\
\hline
\hline
\mathrm{Case} & A_1 & A_2 & A_3 & A_4 & A_5 & A_6 & A_7 & A_8 & A_9 & A_{10} & A_{11} & A_{12} \\
\hline
\end{array}
$$

An example for Case $A_1$ is the following.
Let $\mathcal{C}=\{(1:0:0),(0:1:0),(0:0:1),(1:1:0),(1:0:1)\}$. The two
$4-$secants through $(1:0:0)$ contain the points
$(1:0:0),(0:1:0),(1:1:0),(1:5:0)$ and $(1:0:0),(0:0:1),(1:0:1),(1:0:4),$ respectively.
The $4$-secant through $(0:0:1)$ and $(0:1:0)$ contains the points
$(0:1:1),(0:1:5) \in \mathcal{B}$. The $4$-secant through $(1:1:0)$ and
$(1:0:1)$
contains the points $(1:3:5),(1:2:6) \in \mathcal{B}$.

An example for Case $A_8$ is the following.
Let
$$(1:0:0) \in \mathcal{A} \quad \text{ and }\quad \mathcal{B}=\{(0:1:0),(0:0:1),(0:1:1),(0:1:3)\}.$$
The four $3-$secants
through $(1:0:0)$ contain the points
$$\{(1:0:0),(0:1:0),(1:2:2)\}, \quad \{(1:0:0),(0:1:1),(1:4:4)\},$$
$$\{(1:0:0),(1:1:5),(1:4:6)\}, \quad \{(1:0:0),(1:6:1),(1:3:4)\}.$$

We prove that the other cases do not appear.
Cases $A_5,A_6,A_7A_9,A_{11}$ and $A_{12}$ cannot appear,
because they do not satisfy the condition $a+b+c=11.$
The number of the 4-secants of $\sk $ is $f=(b+2c)/4.$
If $b=0$ and $c=2,$ then $f=1,$ but then obviously do not
exist any point which is on two 4-secants.
In the Cases $A_2,\, A_3$ and $A_4$
we have $f=4.$ But four lines have at most 6 points of intersections,
hence $c=7$ and $c=8$ are impossible.
If $c=6,$ then the four 4-secants form a complete quadrilateral,
the sides of it contain the four points of the set ${\mathcal B},$
and ${\mathcal A}$ consists of a single point, say $P.$
Then each of the four trisecants through $P$ must contain
two points from ${\mathcal B}.$ But then the pigeonhole principle
implies that some of these trisecants have more than one point in common.
This contradiction proves the nonexistence of this configuration.

If $|\sk |=12,$ then Equations (\ref{korlat}) give
$$6x_0+2x_4+6x_5=42 \quad {\rm and} \quad 6x_0+x_2+x_5=48.$$
Proposition \ref{biggerthan12} gives that $x_5\leq 1,$ hence
elementary counting shows that there are only five possibilities
for the numbers $x_0,x_1,\ldots ,x_5.$ These are the following.

\smallskip
\begin{center}
\begin{tabular}{|c|c|c|c|c|c||c|}
\hline
$x_0$&$x_1$&$x_2$&$x_3$&$x_4$&$x_5$&Case\\
\hline
\hline
7&24&6&20&0&0&$B_1$\\
\hline
6&24&12&12&3&0&$B_2$\\
\hline
5&24&18&4&6&0&$B_3$\\
\hline
6&24&11&15&0&1&$B_4$\\
\hline
5&24&17&7&3&1&$B_5$\\
\hline
\end{tabular}
\end{center}

We show that only Case $B_2$ appears.
For each point $P\in \sk $ there are six
secants through $P.$
We have to distribute 11 points among the
secants through $P.$
It follows from Corollary \ref{56} that
$\sk $ has no 6-secant.
Thus the points of $\sk $ can be partitioned into
four subsets. Let ${\mathcal A}\subset \sk $ be the set of points
belonging to five trisecants of $\sk ,$ let ${\mathcal B}\subset \sk $
be the set of points belonging to three trisecants and
one 4-secant of $\sk ,$
let ${\mathcal C}\subset \sk $
be the set of points belonging to one trisecant  and two 4-secants of $\sk $ and finally
let ${\mathcal D}\subset \sk $
be the set of points belonging to two trisecants and one $5$-secant of $\sk .$

If $\sk $ has a $5$-secant, $\ell ,$ then let $\mathcal{R}= \sk \setminus \ell$
and let $\ell \setminus \sk =\{P,Q,R\} .$
\begin{itemize}
\item
Case $B_1$. In this case $\sk$ is a $(12,3)$-arc. In \cite{CS2012} the
intersection sizes with lines of all the regular complete $(12,3)$-arcs in
$\mathrm{PG}(2,7)$ are presented and there exist no regular 
complete $(12,3)$-arcs in
$\mathrm{PG}(2,7)$ having $20$ trisecants. An exhaustive computer  search
among incomplete $(12,3)$-arcs in $\mathrm{PG}(2,7)$ shows that all of them
have less than $20$ 
trisecants.
\item Case $B_3$.
First we prove that no three of the 4-secants have a point in common. There are at most
two 4-secants through any point of $\sk$, and if three 4-secants would meet in a point outside
$\sk ,$ then the union of these lines would contain $\sk ,$ so any other line could contain
at most three points of $\sk ,$ but the total number of $4$-secants is six. Hence through each point
of $\sk$ there are exactly two $4$-secants and one trisecant of $\sk $.

The six 4-secants have
$6\cdot 5/2=15$ points of intersection. Three of these points are not in $\sk ,$ let $X$ and $Y$ be
two of them and let $O$ and $E$ be two points of $\sk $ such that $OX, \, OY, \, EX$ and
$EY$ are 4-secants of $\sk .$ There is a projectivity mapping the points
of the projective frame to $\{ X,Y,O,E\}.$ After this projectivity the points of
$\sk $ are in the affine plane. If we use cartesian
coordinates, we get $O=(0,0), \, E=(1,1),$ and the points $P=OX\cap EY=(1,0)$
and $R=OY\cap EX=(0,1)$ belong to $\sk .$ Let the further points of
$OY\cap \sk $ and $OX\cap \sk $ be $A=(0,a),$ $B=(0,b),$ and $C=(c,0),$ $D=(d,0),$ respectively.
Then $\{ a,b,c,d\} \cap \{ 0,1\} =\emptyset .$

Without loss of generality we may assume that $AC$ and $BD$ are 4-secants.
The equations of these lines
are $X/c+Y/a=1$ and $X/d+Y/b=1,$ respectively.
Then the remaining points of $\sk $ must be $PY\cap AC=K=(1,a-a/c),$
$PY\cap BD=L=(1,b-b/d),$ $RX\cap AC=M=(c-c/a,1)$ and $RX\cap BD=N=(d-d/b,1).$
Hence the lines $OE$ and $PR$ are bisecants. Consider the unique trisecant
through $O.$ It must contain one point from the set $\{ K,L\} $ and one
point from the set $\{ M,N\} .$ But none of the lines $KM$ and $LN$ contains
$O,$ thus without loss of generality we may assume, that the line $KN$ is the
trisecant through $O.$ Hence
\begin{equation}
\label{1}
a-\frac{a}{c}=\frac{1}{d-\frac{d}{b}} \quad \Longleftrightarrow \quad \frac{a(c-1)}{c}=\frac{b}{d(b-1)}.
\end{equation}
In the same way we get that the unique trisecants through the points $P,$ $R$ and $E$
must be the lines $MB,$ $LC$ and $DA,$ respectively. The equation of the line
joining the points $(s,0)$ and $(0,t)$ is $X/s+Y/t=1,$ thus from these collinearity conditions
we get the following equations:
\begin{equation}
\label{2}
\quad \; \; c-\frac{c}{a}+\frac{1}{b}=1 \quad \Longleftrightarrow \quad \frac{c(a-1)}{a}=\frac{b-1}{b},
\end{equation}
\begin{equation}
\label{3}
\quad \; \; \frac{1}{c}+b-\frac{b}{d}=1 \quad \Longleftrightarrow \quad \frac{b(d-1)}{d}=\frac{c-1}{c},
\end{equation}
\begin{equation}
\label{4}
\frac{1}{d}+\frac{1}{a}=1 \quad \Longleftrightarrow \quad d=\frac{a}{a-1}.
\end{equation}

From the last equation we get $(d-1)/d=1/a,$ hence Equations (\ref{3}) and (\ref{1}) give
$ b=(2a-1)/a$ and $bc=1.$ Finally from Equations (\ref{3}) and (\ref{2}) we get $c=(a+1)/a.$ Hence
$$
\frac{2a-1}{a}\cdot \frac{a+1}{a}=1, \quad \mathrm{thus} \quad a^2+a-1=0.$$
But this equation has no root in $\mathrm{GF}(7),$ so there is no semiarc of this type
in $\mathrm{PG}(2,7).$

\item Case $B_4$.
Each point of $\sk \cap \ell$ is contained in two trisecants. Thus the number of trisecants of $\sk$ through the
points of $\ell\cap \sk $ is $10$. Let $x_2^{\prime}$ and  $x_3^{\prime}$ be the number of bisecants and trisecants of
$\mathcal{R}$, respectively. Then counting in two different ways the ordered triples $(A,B, e)$ where
both $A$ and $B$ are points in $\mathcal{R}$ and $e$ is a line incident with both of them, we get
$2x_2^{\prime}+6x_{3}^{\prime}=42.$
On the other hand, each trisecant of $\sk$ containing a point of $\ell$ corresponds to a bisecant of $\mathcal{R}$.
Since the number of trisecants of $\sk$ is $15$, the other $5$ trisecants of $\sk$ must be trisecants also for $\mathcal{R}$, thus $x_2^{\prime}\geq 10$ and $x_{3}^{\prime}=5$. Hence $42=2x_2^{\prime}+6x_{3}^{\prime}\geq 20+30,$ contradiction.
So there is no semiarc of this type.
\item Case $B_5$.
There are no $4$-secants through the points of $\ell \cap \sk .$ Hence each of the three $4$-secants
meets $\mathcal{R}$ in four points. But the union of the three $4$-secants contains at least $4+3+2=9$ distinct points
and $\mathcal{R}$ contains only seven points. So there is no semiarc of this type.
\end{itemize}

Thus only Case $B_2$ can appear.
Now $\sk $ has three 4-secants, say $\ell _1,\ell _2$ and $\ell _3.$
Let $\cal M$ be the set of points of intersections
of the $4$-secants. The number of $4$-secants through any point of $\sk $ is at most two,
hence there are four possibilities.
\begin{enumerate}
\item $|{\cal M}|=1$ and ${\cal M}\cap \sk =\emptyset $,
\item  $|{\cal M}|=3$ and $|{\cal M}\cap \sk |=1 $,
\item $|{\cal M}|=3$ and $|{\cal M}\cap \sk |=2 $,
 \textcolor{red}{}
\item $|{\cal M}|=3$ and $|{\cal M}\cap \sk |=3 $.
\end{enumerate}

An exhaustive computer search shows that there are no examples in
cases 1 and 2, and there are examples in cases
3 and 4. An example of case 3 is the following. Let
$\sk \cap \ell_{1} = \{(1:0:0),(1:0:4),(1:0:5),(1:0:6)\}$,
$\sk \cap \ell_{2} = \{(1:0:0),(0:1:0),(1:5:0),(1:6:0)\}$ and
$\sk \cap \ell_{3} = \{(0:1:0),(0:1:2),(0:1:3),(0:1:5)\}$, finally let
$\sk \setminus (\ell_1 \cup \ell_2 \cup \ell_3)=\{(1:1:1),(1:5:1)\}$.

An example of case 4 is the following.
Let $x=0$, $\sk \cap \ell_{1} = \{(1:0:0),(0:0:1),(1:0:1),(1:0:5)\}$,
$\sk \cap \ell_{2} = \{(1:0:0),(0:1:0),(1:2:0),(1:3:0)\}$ and
$\sk \cap \ell_{3} = \{(0:0:1),(0:1:0),(0:1:4),(0:1:5)\}$, finally let
$\sk \setminus (\ell_1 \cup \ell_2 \cup \ell_3)=\{(1:1:3),(1:1:6),(1:4:3)\}$.
\Qed

Table \ref{table:PG2_7_2} contains the projectively non-equivalent $2$-semiarcs in
$\mathrm{PG}(2,7)$, the number of their $i$-secants, $x_i$,
and the description of the stabilizer groups in $\mathrm{PGL} (3,7)$.

\section{The algorithm}\label{sec:algorithm}


The algorithm used for the classification of $2$-semiarcs in $\mathrm{PG}(2,q)$
is a modification of the one presented in \cite{MMP2007,bar}.
When  possible, the search is helped by the structural constraints proven
in Section \ref{sec:Bounds}.

In this case the algorithm works on \emph{admissible sets}, i.e. sets such
that each point lies on at least two tangent lines, instead of working on
partial solutions. In fact, the property of being a $2$-semiarc is not an
hereditary feature, i.e. a feature conserved by all the subsets, so the weaker
hereditary feature of being an admissible set has been used. It is weaker in
the sense that it allows to prune very few branches of the search space with
respect to the cases when considering arcs and $(k,3)$-arcs. This and the fact
that $2$-semiarcs are in general larger than arcs and $(k,3)$-arcs make the
problem computationally harder than the ones
faced in \cite{MMP2005,MMP2007}.

Note also that, in general, not all the admissible sets can be extended
to $2$-semiarcs.

The exhaustive search has been feasible because projective properties among
admissible sets have been exploited to avoid obtaining too many isomorphic
copies of the same  $2$-semiarc and to avoid searching through parts of the
search space isomorphic to previously searched
ones.

The algorithm starts constructing a tree structure containing a representative
of each class of non-equivalent admissible sets of size less than or equal to
a fixed threshold $h$. If the threshold $h$ were equal to the actual size of
the putative $2$-semiarcs, the algorithm would be orderly, that is capable of
constructing each goal configuration exactly
once \cite{Royle1996}.

However, in the present case, the construction of the tree with the threshold
$h$ equal to the size of the putative $2$-semiarcs would have been too space
and time consuming. For this reason a hybrid approach has been adopted. The
obtained non-equivalent admissible sets of size $h$ have been extended using a
backtracking algorithm trying to determine $2$-semiarcs of the desired
size. In the backtracking phase, the information obtained during the
classification of the admissible sets has been further exploited to prune the
search tree. In fact the points that would have given admissible sets
equivalent to already obtained ones have been excluded from the backtracking
steps.

A simple parallelization technique, based on data distribution, has been used
to divide the load of the computation in a multiprocessor computer. In our
searches we used a $3.3$ Ghz Intel Exacore 16 Gb of
memory.

\section{Results for $8\leq q\leq 13$}\label{sec:8}
\indent
In Table \ref{table:PG2_q<=9}, the number of non-equivalent examples of
$2$-semiarcs in $\mathrm{PG}(2,q)$, $q\leq9$,
is given.
The two examples of $2$-semiarcs of size $8$ in $\mathrm{PG}(2,8)$ are
obtained by deleting two points from the hyperoval (two points of the conic or
one point of the conic and the
nucleus).

The following non-existence results are obvious corollaries of Propositions
\ref{q+1} and \ref{q+2}.
\begin{corollary}\label{cor:qp1}
In $\mathrm{PG}(2,9)$ there are no $2$-semiarcs of size $10$ or $11.$
\end{corollary}

In Tables \ref{table:PG2_8} and \ref{table:PG2_9} the description of the
stabilizer  of the non-equivalent examples of $2$-semiarcs $\sk$ in
$\mathrm{PG}(2,8)$ and $\mathrm{PG}(2,9)$ is presented.   In Table
\ref{table:PG2_8BigStabilizer} (resp. \ref{table:PG2_9BigStabilizer}) the
$2$-semiarcs in $\mathrm{PG}(2,8)$ (resp. $\mathrm{PG}(2,9)$) having
stabilizer of size larger than $16$
are listed
($x_i$ indicates the number of $i$-secants of $\sk$ and $\omega $ denotes an
element satisfying the equation
$\omega^3 + \omega^2 + 1=0$ (resp. $\omega^2 -2\omega - 1=0$)).

By our experimental results we are able to prove the following.
\begin{theorem}
In $\mathrm{PG}(2,11)$ there exist $2$-semiarcs of size $k \in \{11,12,14-26\}$. In $\mathrm{PG}(2,13)$ there exist $2$-semiarcs of size $k \in \{13,27-30\}$.
\end{theorem}
Note that there exists a unique $2$-semiarc of size $11$ (resp. $13$)  in $\mathrm{PG}(2,11)$ (resp. $\mathrm{PG}(2,13)$) and its stabilizer is $(\mathbb{Z}_{11} \rtimes \mathbb{Z}_5) \rtimes \mathbb{Z}_2$ (resp. $(\mathbb{Z}_{13} \rtimes \mathbb{Z}_4) \rtimes \mathbb{Z}_3$), according to Theorem \ref{qcon}. We also proved by an exhaustive computer search that there exists a unique $2$-semiarc of size $12$ in $\mathrm{PG}(2,11)$ and its stabilizer is  $\mathrm{S}_4$.

\bigskip

\bigskip
\noindent
{\bf \Large{Acknowledgement}}

\smallskip
\noindent
The authors are grateful to the anonymous reviewer for his/her
detailed and helpful comments and suggestions.

\begin{center}
\begin{table}
\caption{$2$-semiarcs in $\mathrm{PG}(2,q)$, $q\leq 9$}\label{table:PG2_q<=9}
\begin{center}
\begin{tabular}{|c|c|c|c|}
\hline
$q$&Size&\# non-equivalent examples\\
\hline
\hline
\multirow{3}{*}{$4$}
&$4$&$1$\\
\cline{2-3}
&$6$&$1$\\
\cline{2-3}
&$7$&$1$\\
\hline
\hline
\multirow{3}{*}{$5$}
&$5$&$1$\\
\cline{2-3}
&$6$&$1$\\
\cline{2-3}
&$9$&$1$\\
\hline
\hline
\multirow{4}{*}{$7$}
&$7$&$1$\\
\cline{2-3}
&$9$&$6$\\
\cline{2-3}
&$10$&$12$\\
\cline{2-3}
&$11$&$3$\\
\cline{2-3}
&$12$&$3$\\
\hline
\hline
\multirow{9}{*}{$8$}
&$8$&$2$\\
\cline{2-3}
&$9$&$2$\\
\cline{2-3}
&$10$&$1$\\
\cline{2-3}
&$11$&$10$\\
\cline{2-3}
&$12$&$26$\\
\cline{2-3}
&$13$&$31$\\
\cline{2-3}
&$14$&$29$\\
\cline{2-3}
&$15$&$11$\\
\cline{2-3}
&$16$&$2$\\
\hline
\hline
\multirow{11}{*}{$9$}
&$9$&$1$\\
\cline{2-3}
&$12$&$30$\\

\cline{2-3}
&$13$&$59$\\

\cline{2-3}
&$14$&$360$\\

\cline{2-3}
&$15$&$925$\\

\cline{2-3}
&$16$&$1149$\\
\cline{2-3}
&$17$&$655$\\
\cline{2-3}
&$18$&$162$\\
\cline{2-3}
&$19$&$19$\\
\cline{2-3}
&$20$&$3$\\
\hline
\end{tabular}
\end{center}
\end{table}
\end{center}

\begin{center}
\begin{table}
\caption{$2$-semiarcs in $\mathrm{PG}(2,8)$}\label{table:PG2_8}
\begin{center}
\tabcolsep=0.5 mm
{\scriptsize
\begin{tabular}{|c||c|c|c|c|c|c|c|c|c|c|c|c|c|c|c|c|}
\hline
{\bf Size}&$\mathbb{Z}_{1}$&$\mathbb{Z}_{2}$&$\mathbb{Z}_{3}$&$\mathbb{Z}_{2}^2$&$\mathbb{Z}_{6}$&$\mathcal{S}_{3}$&$\mathbb{Z}_{2}^3$&$\mathbb{Z}_{12}$&$\mathcal{Q}_{6}$&$\mathcal{S}_{3}\times\mathbb{Z}_3$&$\mathrm{S}_4$&$\mathrm{D}_4\times \mathbb{Z}_3$&$\mathrm{A}_4\times \mathbb{Z}_2$&$(\mathbb{Z}_7\rtimes \mathbb{Z}_3)\rtimes \mathbb{Z}_{2}$&$((\mathbb{Z}_4 \times \mathbb{Z}_4) \rtimes \mathbb{Z}_3) \rtimes \mathbb{Z}_2$&$\mathbb{Z}_7 \times (\mathbb{Z}_3 \rtimes \mathbb{Z}_8)$\\
\hline
\hline
{\bf 8}&&&&&&&&&&&&&&$1$&&$1$\\
\hline
{\bf 9}&&&&&&$1$&&&&$1$&&&&&&\\
\hline
{\bf 10}&&&&&&&&&$1$&&&&&&&\\
\hline
{\bf 11}&$5$&$4$&$1$&&&&&&&&&&&&&\\
\hline
{\bf 12}&$8$&$9$&$1$&$1$&&$1$&$1$&$1$&&&$1$&$1$&$2$&&&\\
\hline
{\bf 13}&$22$&$4$&$5$&&&&&&&&&&&&&\\
\hline
{\bf 14}&$14$&$8$&&$6$&&&&&&&&$1$&&&&\\
\hline
{\bf 15}&$5$&$1$&$2$&&$2$&$1$&&&&&&&&&&\\
\hline
{\bf 16}&&&&&&&&&&&&&$1$&&$1$&\\
\hline
\end{tabular}
}
\end{center}
\end{table}
\end{center}

\begin{center}
\begin{table}
\caption{$2$-semiarcs in $\mathrm{PG}(2,9)$}\label{table:PG2_9}
\begin{center}
\tabcolsep=1 mm
{\tiny
\begin{tabular}{|c||c|c|c|c|c|c|c|c|c|c|c|c|c|c|c|c|c|c|c|}
\hline
{\bf Size}&$\mathbb{Z}_{1}$&$\mathbb{Z}_{2}$&$\mathbb{Z}_{3}$&$\mathbb{Z}_{4}$&$\mathbb{Z}_{2}^2$&$\mathbb{Z}_{6}$&$\mathcal{S}_{3}$&$\mathbb{Z}_{2}\times\mathbb{Z}_{4}$&$\mathrm{D}_{4}$&$\mathrm{D}_{6}$&$\mathrm{D}_8$&$\mathrm{D}_4\times \mathbb{Z}_2$&$\mathbb{Z}_4 \rtimes \mathbb{Z}_4$&$\mathcal{S}_{3}\times\mathbb{Z}_3$&$\mathcal{S}_{4}$&$\mathcal{S}_{3}\times\mathbb{Z}_4$&$\mathbb{Z}_6 \times \mathrm{S}_3$&$\mathrm{D}_8 \times \mathrm{S}_3$&$((\mathbb{Z}_3 \times \mathbb{Z}_3) \rtimes \mathbb{Z}_8) \rtimes \mathbb{Z}_2$\\
\hline
\hline
{\bf 9}&&&&&&&&&&&&&&&&&&&$1$\\
\hline
{\bf 12}&$9$&$6$&$1$&&$3$&$4$&$2$&&&$2$&&&&&$1$&$1$&&$1$&\\
\hline
{\bf 13}&$42$&$11$&&$4$&$2$&&&&&&&&&&&&&&\\
\hline
{\bf 14}&$308$&$48$&&$3$&&&&&&&&$1$&&&&&&&\\
\hline
{\bf 15}&$836$&$74$&$3$&&$6$&$2$&$2$&&&$1$&&&&&$1$&&&&\\
\hline
{\bf 16}&$1054$&$73$&&$6$&$11$&&&$2$&$1$&&$1$&&$1$&&&&&&\\
\hline
{\bf 17}&$583$&$59$&&$10$&$1$&&&&&&$1$&&$1$&&&&&&\\
\hline
{\bf 18}&$126$&$22$&$3$&&$3$&&$4$&&&$1$&&&&$1$&$1$&&$1$&&\\
\hline
{\bf 19}&$10$&$5$&&&$2$&&&&$2$&&&&&&&&&&\\
\hline
{\bf 20}&&$2$&&&&&&&&&$1$&&&&&&&&\\
\hline
\end{tabular}
}
\end{center}
\end{table}
\end{center}

\begin{center}
\begin{table}
\caption{$2$-semiarcs in $\mathrm{PG}(2,8)$ with $|G|> 16$}\label{table:PG2_8BigStabilizer}
\begin{center}
\tabcolsep=0.85 mm
\begin{tabular}{|c|c|c|c|c|c|c|c|c|c|}
\hline
$|\mathcal{S}_2|$&$\mathcal{S}_2$&$\ell_{0}$&$\ell_{1}$&$\ell_{2}$&$\ell_{3}$&$\ell_{4}$&$\ell_{5}$&$\ell_{6}$&$G$\\
\hline
\hline
$8$&\tabcolsep=0.85 mm
\begin{tabular}{cccccccc}
 1& 0& 0& 1& 1& 1& 1& 1\\
 0& 1& 0& 1& $\omega$& $\omega^2$& $\omega^3$& $\omega^5$\\
 0& 0& 1& 1& $\omega^5$& $\omega$& $\omega^4$& $\omega^3$\\
\end{tabular}&
29 &
16 &
28 &
0 &
0 &
0 &
0 &
$(\mathbb{Z}_7\rtimes \mathbb{Z}_3)\rtimes \mathbb{Z}_{2}$\\
\hline
$8$&
\tabcolsep=0.85 mm
\begin{tabular}{cccccccc}
 1& 0& 0& 1& 1& 1& 1& 1\\
 0& 1& 0& 1& $\omega$& $\omega^3$& $\omega^5$& $\omega^6$\\
 0& 0& 1& 1& $\omega^5$& $\omega^6$& $\omega^4$& $\omega$\\
\end{tabular}&
29 &
16 &
28 &
0 &
0 &
0 &
0 &
$\mathbb{Z}_7 \times (\mathbb{Z}_3 \rtimes \mathbb{Z}_8)$\\
\hline
\hline
$9$&\tabcolsep=0.85 mm
\begin{tabular}{ccccccccc}
 1& 0& 0& 1& 0& 1& 1& 1& 1\\
 0& 1& 0& 1& 1& $\omega$& $\omega^2$& $\omega^3$& $\omega^6$\\
 0& 0& 1& 1& $\omega$& $\omega^5$& $\omega^4$& $\omega^6$& $\omega^2$\\
\end{tabular}&
25 &
18 &
27 &
3 &
0 &
0 &
0 &
$\mathrm{S}_{3}\times \mathbb{Z}_{3}$\\
\hline
\hline
$12$&\tabcolsep=0.85 mm
\begin{tabular}{cccccccccccc}
 1& 0& 0& 1& 0& 0& 0& 0& 1& 1& 1& 1\\
 0& 1& 0& 1& 1& 1& 1& 1& $\omega$& $\omega^2$& $\omega^3$& $\omega^5$\\
 0& 0& 1& 1& $\omega$& $\omega^2$& $\omega^3$& $\omega^5$& $\omega$& $\omega^2$& $\omega^3$& $\omega^5$\\
\end{tabular}&
11 &
24 &
36 &
0 &
0 &
0 &
2 &
$\mathrm{S}_{4}$\\
\hline

$12$&\tabcolsep=0.85 mm
\begin{tabular}{cccccccccccc}
 1& 0& 0& 1& 0& 0& 1& 1& 1& 1& 1& 1\\
 0& 1& 0& 1& 1& 1& 0& 0& $\omega$& $\omega^2$& $\omega^4$& $\omega^5$\\
 0& 0& 1& 1& $\omega$& $\omega^5$& $\omega$& $\omega^2$& $\omega^5$& $\omega^5$& $\omega^5$& $\omega^6$\\
\end{tabular}&
13 &
24 &
30 &
0 &
6 &
0 &
0 &
$\mathrm{D}_{4}\times \mathbb{Z}_{3}$\\
\hline

$12$&\tabcolsep=0.85 mm
\begin{tabular}{cccccccccccc}
 1& 0& 0& 1& 0& 0& 1& 1& 1& 1& 1& 1\\
 0& 1& 0& 1& 1& 1& 0& $\omega$& $\omega^2$& $\omega^4$& $\omega^5$& $\omega^5$\\
 0& 0& 1& 1& $\omega$& $\omega^5$& $\omega^2$& $\omega^5$& $\omega^5$& 1& $\omega$& $\omega^6$\\
\end{tabular}&
14 &
24 &
24 &
8 &
3 &
0 &
0 &
$\mathrm{A}_{4}\times\mathbb{Z}_2$\\
\hline

$12$&\tabcolsep=0.85 mm
\begin{tabular}{cccccccccccc}
 1& 0& 0& 1& 0& 0& 1& 1& 1& 1& 1& 1\\
 0& 1& 0& 1& 1& 1& 0& $\omega$& $\omega$& $\omega^2$& $\omega^5$& $\omega^6$\\
 0& 0& 1& 1& $\omega$& $\omega^5$& $\omega^2$& $\omega^4$& $\omega^5$& $\omega^5$& $\omega^6$& $\omega^2$\\
\end{tabular}&
14 &
24 &
24 &
8 &
3 &
0 &
0 &
$\mathrm{A}_{4}\times\mathbb{Z}_2$\\
\hline
\hline
$14$&\tabcolsep=0.85 mm
\begin{tabular}{cccccccccccccc}
 1& 0& 0& 1& 0& 0& 0& 0& 1& 1& 1& 1& 1& 1\\
 0& 1& 0& 1& 1& 1& 1& 1& 0& $\omega^2$& $\omega^3$& $\omega^3$& $\omega^5$& $\omega^6$\\
 0& 0& 1& 1& 1& $\omega$& $\omega^2$& $\omega^5$& $\omega^2$& $\omega^4$& $\omega$& $\omega^6$& $\omega^6$& 0\\
\end{tabular}&
6 &
28 &
25 &
12 &
0 &
0 &
2 &
$\mathrm{D}_{4}\times\mathbb{Z}_3$\\
\hline
\hline
$16$&\tabcolsep=0.85 mm
\begin{tabular}{cccccccccccccccc}
 1& 0& 0& 1& 0& 1& 1& 1& 1& 1& 1& 1& 1& 1& 1& 1\\
 0& 1& 0& 1& 1& 0& 1& $\omega$& $\omega$& $\omega$& $\omega^2$& $\omega^2$& $\omega^4$& $\omega^4$& $\omega^5$& $\omega^5$\\
 0& 0& 1& 1& $\omega$& $\omega^2$& $\omega^2$& 0& 1& $\omega^4$& $\omega^2$& $\omega^5$& $\omega^5$& $\omega^6$& $\omega^4$&$\omega^6$\\
\end{tabular}&
5 &
32 &
0 &
32 &
4 &
0 &
0 &
$\mathrm{A}_{4}\times\mathbb{Z}_2$\\
\hline

$16$&\tabcolsep=0.85 mm
\begin{tabular}{cccccccccccccccc}
 1& 0& 0& 1& 0& 1& 1& 1& 1& 1& 1& 1& 1& 1& 1& 1\\
 0& 1& 0& 1& 1& 0& 0& 1& $\omega$& $\omega$& $\omega^2$& $\omega^2$& $\omega^3$& $\omega^3$& $\omega^5$& $\omega^5$\\
 0& 0& 1& 1& $\omega$& $\omega^2$& $\omega^5$& $\omega^4$& 0& 1& $\omega^2$& $\omega^5$& 1& $\omega^6$& $\omega^4$&$\omega^6$\\
\end{tabular}&
5 &
32 &
0 &
32 &
4 &
0 &
0 &
$((\mathbb{Z}_4 \times \mathbb{Z}_4) \rtimes \mathbb{Z}_3) \rtimes \mathbb{Z}_2$\\
\hline
\end{tabular}
\end{center}
\end{table}
\end{center}

\begin{center}
\begin{table}
\caption{$2$-semiarcs in $\mathrm{PG}(2,9)$ with $|G|> 16$}\label{table:PG2_9BigStabilizer}
\begin{center}
\tabcolsep=0.85 mm
\begin{tabular}{|c|c|c|c|c|c|c|c|c|}
\hline
$|\mathcal{S}_2|$&$\mathcal{S}_2$&$\ell_{0}$&$\ell_{1}$&$\ell_{2}$&$\ell_{3}$&$\ell_{4}$&$\ell_{5}$&$G$\\
\hline
\hline
$9$&\tabcolsep=0.85 mm
\begin{tabular}{ccccccccc}
 1& 0& 0& 1& 1& 1& 1& 1& 1\\
 0& 1& 0& 1& $\omega$& 2& $\omega^5$& $\omega^6$& $\omega^7$\\
 0& 0& 1& 1& $\omega^5$& $\omega$& $\omega^3$& $\omega^7$& $\omega^2$\\
\end{tabular}&
37 &
18 &
36 &
0 &
0 &
0&
$((\mathbb{Z}_3 \times \mathbb{Z}_3) \rtimes \mathbb{Z}_8) \rtimes \mathbb{Z}_2$\\
\hline
\hline
$12$&
\tabcolsep=0.85 mm
\begin{tabular}{cccccccccccc}
1& 0& 0& 1& 0& 0& 1& 1& 1& 1& 1& 1\\
 0& 1& 0& 1& 1& 1& 1& $\omega$& $\omega^3$& 2& $\omega^5$& $\omega^6$\\
 0& 0& 1& 1& $\omega$& $\omega^2$& $\omega^3$& $\omega^7$& $\omega^6$& $\omega^5$& $\omega$& $\omega$\\
\end{tabular}&
24 &
24 &
36 &
4 &
3 &
0&
$\mathrm{S}_{3}\times\mathbb{Z}_4$\\
\hline
$12$&\tabcolsep=0.85 mm
\begin{tabular}{cccccccccccc}
 1& 0& 0& 1& 0& 1& 1& 1& 1& 1& 1& 1\\
 0& 1& 0& 1& 1& 0& $\omega$& $\omega^2$& $\omega^3$& 2& 2& $\omega^6$\\
 0& 0& 1& 1& $\omega$& 2& $\omega^7$& $\omega^7$& $\omega^3$& 1& $\omega^2$& $\omega$\\
\end{tabular}&
25 &
24 &
30 &
12 &
0 &
0&
$\mathrm{S}_{4}$\\
\hline
$12$&\tabcolsep=0.85 mm
\begin{tabular}{cccccccccccc}
 1& 0& 0& 1& 0& 0& 1& 1& 1& 1& 1& 1\\
 0& 1& 0& 1& 1& 1& 0& $\omega$& $\omega^3$& 2& $\omega^5$& $\omega^7$\\
 0& 0& 1& 1& $\omega$& $\omega^5$& $\omega^3$& $\omega^7$& $\omega$& $\omega^2$& $\omega^7$& $\omega^6$\\
\end{tabular}&
24 &
24 &
36 &
4 &
3 &
0&
$\mathrm{D}_8 \times \mathrm{S}_3$\\
\hline
\hline
$15$&\tabcolsep=0.85 mm
\begin{tabular}{ccccccccccccccc}
 1& 0& 0& 1& 0& 1& 1& 1& 1& 1& 1& 1& 1& 1& 1\\
 0& 1& 0& 1& 1& 0& 1& $\omega$& $\omega$& $\omega^3$& $\omega^3$& $\omega^5$& $\omega^5$& $\omega^6$& $\omega^7$\\
 0& 0& 1& 1& $\omega$& $\omega^5$& $\omega^6$& 2& $\omega^7$& 1& $\omega^6$& $\omega$& $\omega^5$& $\omega^7$& 2\\
\end{tabular}&
16 &
30 &
15 &
30 &
0 &
0&
$\mathrm{S}_{4}$\\
\hline
\hline
$18$&\tabcolsep=0.85 mm
\begin{tabular}{cccccccccccccccccc}
 1& 0& 0& 1& 0& 0& 0& 1& 1& 1& 1& 1& 1& 1& 1& 1& 1& 1\\
 0& 1& 0& 1& 1& 1& 1& 0& $\omega$& $\omega$& $\omega$& $\omega$& $\omega^3$& $\omega^3$& $\omega^5$& $\omega^6$& $\omega^7$& $\omega^7$\\
 0& 0& 1& 1& 1& $\omega$& $\omega^3$& $\omega^3$& 1& $\omega$& 2& $\omega^7$& $\omega^2$& $\omega^3$& 1& 1& 0&  $\omega^2$\\
\end{tabular}&
1 &
36 &
36 &
9 &
0 &
9 &
$\mathrm{S}_{3}\times\mathbb{Z}_3$\\
\hline
$18$&\tabcolsep=0.85 mm
\begin{tabular}{cccccccccccccccccc}
 1& 0& 0& 1& 0& 0& 1& 1& 1& 1& 1& 1& 1& 1& 1& 1& 1& 1\\
 0& 1& 0& 1& 1& 1& 0& 1& $\omega$& $\omega$& $\omega$& $\omega^3$& $\omega^3$& $\omega^3$& 2& $\omega^6$& $\omega^6$& $\omega^7$\\
 0& 0& 1& 1& $\omega$& $\omega^5$& 1& $\omega^3$& 2& $\omega^5$& $\omega^7$& 0& $\omega^3$& $\omega^7$& 1& $\omega$& 2& 2\\
\end{tabular}&
6 &
36 &
15 &
22 &
12 &
0&
$\mathrm{S}_{4}$\\
\hline
$18$&\tabcolsep=0.85 mm
\begin{tabular}{cccccccccccccccccc}
 1& 0& 0& 1& 0& 1& 1& 1& 1& 1& 1& 1& 1& 1& 1& 1& 1& 1\\
 0& 1& 0& 1& 1& 0& 0& 1& 1& $\omega$& $\omega$& $\omega$& $\omega^2$& $\omega^3$& 2& 2& 2& $\omega^5$\\
 0& 0& 1& 1& $\omega$& 1& $\omega^7$& $\omega^2$& $\omega^3$& 2& $\omega^6$& $\omega^7$& 1& 0& 0& $\omega^2$& $\omega^7$& $\omega^2$\\
\end{tabular}&
4 &
36 &
27 &
6 &
18 &
0&
$\mathbb{Z}_6 \times \mathrm{S}_3$\\
\hline
\end{tabular}
\end{center}
\end{table}
\end{center}

\begin{flushleft}
Gy\"orgy Kiss \\
Department of Geometry and
MTA-ELTE GAC Research Group \\
E\"otv\"os Lor\'and University \\
1117 Budapest, P\'azm\'any s. 1/c, Hungary \\
e-mail: {\sf kissgy@cs.elte.hu}
\end{flushleft}

\begin{flushleft}

Daniele Bartoli, Giorgio Faina, Stefano Marcugini and Fernanda Pambianco \\
Dipartimento di Matematica e Informatica,
Universit\`{a} degli Studi di Perugia \\
Via Vanvitelli 1, 06123 Perugia, Italy \\
e-mails: {\sf daniele.bartoli@dmi.unipg.it}, {\sf faina@dmi.unipg.it},
{\sf gino@dmi.unipg.it},
{\sf fernanda@dmi.unipg.it}
\end{flushleft}
\end{document}